\documentclass[reqno,11pt]{amsart}

\usepackage{amssymb}

\makeatletter

\newtheoremstyle{sftheorem}
  {5pt}
  {5pt}
  {\itshape}
  {}
  {\sffamily\mdseries\upshape}
  {.}
  {.5em}
  {}

 \theoremstyle{sftheorem}

 \newtheorem{theorem}{Theorem}

 \newtheorem{lemma}{Lemma}[section]
 \newtheorem{proposition}[lemma]{Proposition}

 \newcommand{\df}[1]{{\rmfamily\itshape\mdseries#1}}
 \newcommand\setof[1]{\mathopen\{\,#1\,\mathclose\}}
 \newcommand\Prob{\mathop{\operator@font Prob}}
 \newcommand\cei[1]{{\lceil #1\rceil}}
 
 \newcommand\flo[1]{{\lfloor #1\rfloor}}
 \let\et=\wedge

 \newcommand{\ZZ}{\mathbb{Z}}
 \newcommand \ag{\alpha}
 \newcommand{\cE}{\mathcal{E}}
 \newcommand{\cF}{\mathcal{F}}
 \newcommand{\cG}{\mathcal{G}}

\makeatother

\begin{document}

	\title{The Clairvoyant Demon Has a Hard Task}

\author{Peter G\'acs}
\address{Computer Science Department\\ Boston University}
\email{gacs@bu.edu}

\begin{abstract}
  Consider the integer lattice $L = \ZZ^{2}$.
  For some $m\ge 4$, let us color each column of this lattice independently
and uniformly into one of $m$ colors.
  We do the same for the rows, independently from the columns.
  A point of $L$ will be called \df{blocked} if its row and column have the
same color.
  We say that this random configuration \df{percolates} if there is a path
in $L$ starting at the origin, consisting of rightward and upward unit
steps, and avoiding the blocked points.
  As a problem arising in distributed computing, it has been conjectured
that for $m\ge 4$, the configuration percolates with positive probability.
  This has now been proved (in a later paper) for large $m$.
We prove that the probability that there
is percolation to distance $n$ but not to infinity is not exponentially
small in $n$.
  This narrows the range of methods available for proving the conjecture.
\end{abstract}

\keywords{Dependent percolation, scheduling, distributed computing}

\maketitle

  \section{Statement of the result}

\subsection{Introduction}

  Let $x = (x(0),x(1),\ldots)$ be an infinite sequence and
$u=(u(0),u(1),\ldots)$ be a binary sequence with elements in
$\{0,1\}$.
  Let $s_{n}=\sum_{i=0}^{n-1}u(i)$.
  We define the \df{delayed version} $x^{(u)}$ of $x$, by
 \[
   x^{(u)}(n) = x(s_{n}).
 \]
  Thus, if $i=s_{n}$ then $x^{(u)}(n)=x(i)$, and $x^{(u)}(n+1)=x(i)$ or
$x(i+1)$ depending on whether $u(n)=0$ or 1.
  If $u(n)=0$ then we can say that $x^{(u)}$ is \df{delayed} at time $n+1$.
  For two infinite sequences $x$, $y$ we say that they \df{do not collide}
if there is a delay sequence $u$ such that for each $n$ we have
 \[
   x^{(u)}(n)\ne y^{(1-u)}(n).
 \]
  Here, $1-u$ is the delay sequence complementary to $u$: thus, $y^{(1-u)}$
is delayed at time $n$ if and only if $x^{(u)}$ is not.

  For a given $m>1$, suppose that $X=(X(0),X(1),\ldots)$ is an
infinite sequence of independent random variables, and
 $Y = (Y(0),Y(1),\ldots)$ is another such sequence, also independent of
$X$, where all variables are uniformly distributed over $\{0,\ldots,m-1\}$.
The following theorem has been conjectured in~\cite{CopperTetWink93}.

  \begin{proposition}[See \cite{GacsWalks01}]
  If $m$ is sufficiently large then with positive probability, $X$ does not
collide with $Y$.
  \end{proposition}

  The sequences $X,Y$ can be viewed as two independent random walks on the
complete graph $K_{m}$, and then the problem is whether a ``clairvoyant
demon'', i.e.~ a being who knows in advance both infinite sequences $X$ and
$Y$, can introduce delays into these walks in such a way that they never
collide.

\subsection{Graph reformulation}

  We define a graph $G = (V,E)$ as follows.
  $V=\ZZ_{>0}^{2}$ is the set of points $(i,j)$ where $i,j$ are positive
integers.
  Let us define the \df{distance} of two points $(i,j)$, $(k,l)$ as
$|k-i|+|l-j|$ ($L_{1}$ distance).
  When representing the set $V$ of points $(i,j)$ graphically, the right
direction is the one of growing $i$, and the upward direction is the one of
growing $j$.
  The set $E$ of edges consists of all pairs of the form $((i,j),(i+1,j))$
and $((i,j),(i,j+1))$.

  Given $X,Y$ as in the theorem, let us say that a point $(i,j)$ has color
$k$ if $X(i)=Y(j)=k$.
  Otherwise, it has color $-1$, which we will call \df{white}.
  It is easy to see that $X$ and $Y$ do not collide if and only if there is
an infinite directed path in $G$ starting from $(0, 0)$ and proceeding on
white points.
  Indeed, each path corresponds to a delay sequence $u$ such that $u(n)=1$
if and only if the edge is horizontal.
  Thus, the two sequences do not collide if and only if the graph of white
points ``percolates''.
  We will say that \df{there is percolation} if the probability that there
is an infinite path is positive.

  It has been shown independently in \cite{BalBollobStacCperc99} and
\cite{WinklerCperc99} that if the graph is undirected then there is
percolation even for $m = 4$.
  In traditional percolation theory, when there is percolation then
typically (unless the probability of blocking is at a ``critical
point'') the probability that there is percolation to a distance $n$ but no
percolation to infinity is exponentially small in $n$.
  This is the case also in the papers cited above, but it is not true for
the directed percolation problem we are facing.

  \begin{theorem}\label{t.lb}
  If there is percolation from the origin to infinity with positive
probability then the probability of percolating from the origin to distance
$n$ but not to infinity is at least $Cn^{-\ag}$ for some constants
$C,\ag > 0$ depending only on $m$.
 \end{theorem}

\section{The proof}

  Let $b_{m}=(0,1,2,\ldots,m-1)$ be called the \df{basic color sequence} of
length $m$: it is simply the list of all different colors.
  Let $b'_{m}$ be the reverse of $b_{m}$, i.e.~$b'_{m}(i)=b_{m}(m-i-1)$.
  Let $\cE_{n, k}$ be the event that for all $i\in [0, k - 1]$,
$j \in [0, m - 1]$ we have
 \[
   Y(n + i m + j - 1) = b'_{m}(j), 
 \]
  i.e.~starting with the index $n - 1$, the sequence $Y$ has $k$
consecutive repetitions of $b'_{m}$.
  We say that $i$ is an index of the occurrence of $b_{m}$ in the sequence
$X$ if $X(i + j)=b_{m}(j)$ for $j\in[0, m - 1]$.
  For $i > 0$, let $\tau_{i}$ be the $i$-th index of occurrence of $b_{m}$
in $X$.
  Let $\cF_{n, k}$ be the event that for all $i \in [1, n]$ we have
 \begin{equation}\label{e.taucond}
  \tau_{i + 1} - m - \tau_{i} \le k - 1,
 \end{equation}
 and also $\tau_{1} \le k - 1$.

  \begin{lemma}
  If for some integer $n>0$ both $\cE_{n,k}$ and $\cF_{n,k}$ hold then
there is no white infinite directed path.
  \end{lemma}

  \begin{proof}
  Let us assume that there is a white infinite path.
  Since there are $n$ consecutive copies of $b_{m}$ in $X$, the $i$-th one
starting at index $\tau_{i}$, for each $p\in[1,n]$ there must be a vertical
step in the path with an $x$ projection in $[\tau_{p}, \tau_{p} + m - 1]$.
  Therefore the path ascends to the segment
 $[0, \tau_{n} + m - 1] \times \{n - 1\}$, before its $x$ projection
reaches $t_{n}$. 

  For each $1 \le p\le n$, $0 \le q < k$ there is a diagonally descending
sequence of $m$ colored points
 \[
  \setof{(\tau_{p} + j,\; n + (q + 1)m - j - 2) : j \in [0, m - 1]}.
 \]
  For a fixed $p$ there are $k$ such diagonal barriers stacked above each
other, forming an impenetrable column of height $k m$.
  The path would have to ascend between two of these columns, say between
column $i$ and $i + 1$.
  The distance of two consecutive columns from each other is 
 \[
 \tau_{i+1} - m - \tau_{i} \le k - 1.
 \]
  Since there are $k$ consecutive copies of $b'_{m}$ in $Y$ starting at
index $n - 1$, for each $q\in[0, k - 1]$ there must be a horizontal step in
the path with a height in $n + qm - 1 + [0, m - 1]$.
  Since the distance of the two columns is at most $k - 1$, it is not
possible for the path to pass between the two columns.
  (The same holds for the space before the first column.)
  \end{proof}

 \begin{lemma}\label{l.lb.prob-bounds}
  There is a constant $\ag$ such that for all $s > 0$, there is a $k$ with
  \[
    \Prob(\cE_{n,k}) \ge m^{-m}n^{-\ag (s+1)},
\quad  \Prob(\cF_{n,k}) \ge 1 - n^{-s}.
  \]
 \end{lemma}

 \begin{proof}  With $p_{1} = m^{-m}$, we have
 $\Prob(\cE_{n,k}) =  p_{1}^{k}$.
  Let us estimate the probability of $\cF_{n, k}$.
  The probability that $\tau_{1} > k - 1$ is upperbounded by the
probability that a copy of $b_{m}$ does not begin at $i$ in $X$ for any $i$
in $\{0, m, \ldots, \flo{(k - 1) / m} m\}$, which is
 \[
   (1 - p_{1})^{\flo{(k - 1) / m} + 1} \le (1 - p_{1})^{k / m} < 
	e^{-p_{1} k / m}.
 \]
  The same estimate holds for the probability of ~\eqref{e.taucond}
assuming that the similar conditions for smaller $i$ have already been
satisfied.
  Hence
 \[
  \Prob(\cF_{n, k}) > 1 - n e^{-p_{1}k/m}.
 \]
  Let us choose
 \begin{equation}\label{e.s2k}
   k = \cei{(s+1) (m \log n)/p_{1}}
 \end{equation}
  for some $s > 0$, then $\Prob(\cF_{n,k}) > 1 - n^{-s}$, while
 \[
  \Prob(\cE_{n,k}) \ge p_{1}n^{\frac{m \log p_{1}}{p_{1}}(s + 1)}.
 \]
 \end{proof}

 \begin{proof}[Proof of Theorem \protect\ref{t.lb}]
  Assume that there is an infinite path with some positive probability
$P_{0}$.
  Choose $s$ such that $n^{-s}< 0.5 P_{0}$ and choose $k$ as a function of
$s$ as in ~\eqref{e.s2k}.
  Then the probability that $\cF_{n, k}$ holds and there is an infinite path
is at least $0.5 P_{0}$.
  Let $\cG_{n}$ be the event that there is a path leaving the rectangle
$[0, \tau_{n}]\times[0, n-1]$.
  Then $\Prob(\cF_{n, k}\et \cG_{n})\ge 0.5 P_{0}$.
  Since $\cF_{n, k}\et \cG_{n}$ is independent of $\cE_{n,k}$, we have
 \[
   \Prob(\cE_{n,k}\et \cF_{n, k}\et \cG_{n})
	 \ge 0.5 P_{0}m^{-m}n^{-\ag (s+1)}.
 \]
 \end{proof}

  The author is thankful to J\'anos Koml\'os and Endre Szemer\'edi for
telling him about the problem.

\end{document}